\documentclass[journal]{IEEEtran}
\usepackage{amsmath}
\usepackage{graphicx}
\usepackage{dsfont}
\usepackage{epstopdf}
\usepackage{subfigure}
\usepackage{algorithm}
\usepackage{algpseudocode}
\usepackage{cite}
\usepackage{multicol}
\usepackage{multirow}
\usepackage{mathtools}
\usepackage{amssymb}
\usepackage{color}
\usepackage{bm}
\usepackage{indentfirst}
\usepackage{url}
\usepackage{nomencl}
\usepackage{booktabs}
\usepackage{cite}
\usepackage[numbers,sort&compress]{natbib}
\usepackage{graphicx}
\usepackage{subfigure}
\usepackage{algorithm}
\usepackage{algpseudocode}
\usepackage{cite}
\usepackage{multicol}
\usepackage{multirow}
\usepackage{mathtools}
\usepackage{amssymb}
\usepackage{color}
\usepackage{bm}
\usepackage{indentfirst}
\usepackage{url}
\usepackage{nomencl}
\usepackage{booktabs}
\usepackage{cite}
\usepackage[numbers,sort&compress]{natbib}
\usepackage{amsthm}
\newtheorem{Theorem}{Theorem}
\newtheorem{Lemma}{Lemma}

\newtheorem{Remark}{Remark}

\newtheorem{Assumption}{Assumption}
\allowdisplaybreaks

\begin{document}

\title{Distributed Optimization with Coupling Constraints Based on Dual Proximal Gradient Method in Multi-Agent Networks} 


\author{Jianzheng Wang, Guoqiang~Hu,~\IEEEmembership{Senior Member,~IEEE} 
	\thanks{This work was supported in part by Singapore Economic Development Board under EIRP grant S14-1172-NRF EIRP-IHL, and in part by the Republic of Singapore s National Research Foundation under its Campus for Research Excellence and Technological Enterprise (CREATE) Programme through a grant to the Berkeley Education Alliance for Research in Singapore (BEARS) for the Singapore-Berkeley Building Efficiency and Sustainability in the Tropics (SinBerBEST) Program.}
	\thanks{Jianzheng Wang and Guoqiang Hu are with the School of Electrical and Electronic Engineering, Nanyang Technological University, Singapore, 639798 e-mail: (wang1151@e.ntu.edu.sg, gqhu@ntu.edu.sg).}
}
\maketitle

\begin{abstract}                          
In this paper, we aim to solve a distributed optimization problem with affine coupling constraints in a multi-agent network, where the cost function of the agents is composed of smooth and possibly non-smooth parts. To solve this problem, we resort to the dual problem by deriving the Fenchel conjugate, resulting in a consensus-based constrained optimization problem. Then, we propose a distributed dual proximal gradient algorithm, where the agents make decisions based on the information of immediate neighbors. Provided that the non-smooth parts in the primal cost functions are with some simple structures, we only need to update dual variables by some simple operations, by which the overall computational complexity can be reduced. An ergodic convergence rate of the proposed algorithm is derived and the feasibility is numerically verified by solving a social welfare optimization problem in the electricity market.
\end{abstract}


\section{Introduction}

\subsection{Background and Motivation}

\IEEEPARstart{D}{ecentralized} optimization has become an active topic in recent years for solving various engineering problems, such as detection and localization in sensor networks \cite{lesser2003distributed}, machine learning problems \cite{lee2017speeding}, and economic dispatch in power systems \cite{bai2017distributed}. As a typical optimization procedure, each agent usually maintains an individual decision variable and the global optimal solution can be obtained with multiple rounds of communication and decision making. In this paper, we focus on a class of composite optimization problems, where the cost functions are composed of smooth (differentiable) and possibly non-smooth (non-differentiable) parts, which are often discussed in resource allocation problems \cite{beck2014fast}, least absolute shrinkage and selection operator (LASSO) regressions \cite{hans2009bayesian}, and support vector machines \cite{zhao2017scope}. To solve these problems, widely discussed techniques include alternating direction method of multipliers, primal-dual subgradient method, and proximal gradient method, etc.

The majority of existing works on decentralized optimizations assume that the agents are fully connected to ensure the correctness of the optimization results, which limits their usage in large-scale distributed networks \cite{beck2009fast}. To overcome this issue, a valid alternative is applying graph theory in modeling the communication links, leading to the distributed setup where the agents only communicate with their immediate neighbors \cite{yang2019survey}. However, with the increasing demand in the computational efficiency of various fields, more explorations on the algorithm development for distributed optimization problems (DOPs) are required \cite{bertsekas1989parallel}. Observing that proximal gradient method is usually numerically more stable than the subgradient-based counterpart in composite optimization problems \cite{bertsekas2011incremental}, in this work, we aim to develop an efficient distributed optimization algorithm based on proximal gradient method.


\subsection{Literature Review}

Fruitful distributed algorithms for solving DOPs can be found in the existing works. To adapt to large-scale distributed networks, consensus-based DOPs without coupling constraint were studied in \cite{pang2019randomized,nedic2009distributed}, where the agents make decisions with local variables and certain agreement on the optimal solution is achieved through local communication.
Alternatively, we focus on optimizing a class of composite DOPs subject to affine coupling constraints, where the global objective functions are the sum of the local cost functions.
To solve the problems of interest, some primal-dual subgradient methods were studied in the past few years, where the average consensus technique was employed for the iteration of local variables \cite{zhu2011distributed,chang2014distributed,yuan2011distributed}.
Alternatively, by formulating some local subproblems in the iteration steps, dual decomposition methods were actively investigated in the recent works, where the agents exchange the dual information with their neighbors based on the result of subproblems \cite{necoara2017fully,notarnicola2019constraint,falsone2016distributed,falsone2017dual,simonetto2016primal,necoara2015linear}. When it comes to composite optimizations, dual proximal gradient method was exploited recently by applying the proximal gradient method to the dual problems as in \cite{notarnicola2016asynchronous,beck2014fast,kim2016fast,wang2021distributed}, where, however, no general affine coupling constraints were considered in the distributed computations.

Different from the aforementioned works, we aim to incorporate the dual proximal gradient method in a distributed setup with general affine coupling constraint. Then, a distributed dual proximal gradient (DDPG) algorithm is proposed, by which the computational complexity of the local computations can be reduced if the non-smooth parts of the cost functions are with some simple structures. To highlight of features of the DDPG algorithm, the comparisons with some state-of-the-art works with similar problem setups are summarized as follows.
\begin{itemize}
\item  One key feature of the proposed DDPG algorithm is that one only needs to update the dual variables by some simple operations provided that the proximal mapping of the non-smooth parts in the primal cost functions can be explicitly derived, which is different from \cite{necoara2017fully,notarnicola2019constraint,su2021distributed,li2020distributed,falsone2016distributed,falsone2017dual,chang2016proximal,simonetto2016primal}.\footnote[1]{For the DOPs with smooth cost functions, some existing works on dual algorithms, e.g., \cite{necoara2015linear}, also can avoid the update of primal variables. However, directly extending their results to non-smooth cases can be costly in the sense that the computation of the gradient of the formulated dual functions requires an additional nontrivial optimization process. Therefore, the contribution to the computational efficiency of this work is established for possibly non-smooth cost functions.} This settlement can be efficient in the sense that no costly inner-loop optimizations of the primal or other auxiliary variables are required. In addition, compared with the algorithms without inner-loop optimizations, no explicit convergence rates were provided in \cite{zhu2011distributed,chang2014distributed}. By contrast, an asymptotic convergence is ensured by the DDPG algorithm with an ergodic convergence rate $\mathcal{O}(1/T)$ for the dual function value.
\item 
        In terms of the mathematical assumptions, the algorithms in \cite{zhu2011distributed,li2020distributed,su2021distributed,notarnicola2019constraint,falsone2017dual,chang2014distributed,simonetto2016primal,falsone2016distributed} require some compact constraints on the primal variables and \cite{yuan2011distributed,mateos2016distributed} assume some bounded subgradients to ensure the convergence of their algorithms. By contrast, this work focuses on dual sequences without the boundedness requirement on the primal variables or subgradients. Also, the objective functions in \cite{chang2014distributed,zhu2011distributed} are assumed to be continuous, while the DDPG algorithm allows the non-smooth parts of the objective functions to be lower semi-continuous.
\end{itemize}

The contributions of this work are summarized as follows.
\begin{itemize}
  \item We consider a class of composite DOPs with affine coupling constraints. A DDPG algorithm is proposed by deriving the dual problem based on Fenchel conjugate, where the optimal solution can be obtained when the agents update only with the dual information of immediate neighbors, leading to a distributed computation environment.
  \item The proposed DDPG algorithm only requires the update of dual variables by some simple operations if the non-smooth parts of the cost functions are simple-structured, which can reduce the overall computational complexity. Then, some cost functions with more complicated structures are discussed by tolerating some inner-loop optimizations in each iteration. An ergodic convergence rate of the DDPG algorithm is derived.
\end{itemize}

The remainder of this paper is organized as follows. Section \ref{sa2} provides some fundamental definitions and mathematical properties employed by this work. Section \ref{sa3} formulates the optimization problem of interest and introduces the assumptions. In Section \ref{sa4}, the DDPG algorithm is proposed based on the dual problem. The convergence analysis is conducted in Section \ref{sa5}. The feasibility of the proposed algorithm is verified by a numerical simulation in Section \ref{sa6}. Section \ref{sa7} concludes this paper.

\section{Preliminaries}\label{sa2}

Some frequently used notations, definitions, and relevant properties of graph theory, proximal mapping, and Fenchel conjugate are provided in this section.

\subsection{Notations}
$\mathbb{N}$ and $\mathbb{N}_+$ denote the non-negative and positive integer spaces, respectively. Let $\mid {A}\mid$ be the size of set ${A}$.
Operator $(\cdot)^{\top}$ represents the transpose of a matrix.
${A}_1 \times {A}_2$ denotes the Cartesian product of sets ${A}_1$ and ${A}_2$.
$\mathbf{relint}{A}$ represents the relative interior of set ${A}$. 
Let $\| \mathbf{u}\|^2_{\mathbf{X}}= \mathbf{u}^{\top} \mathbf{X}\mathbf{u}$ with $\mathbf{X}$ being a square matrix. $\mathbf{X} \succ 0$ ($\succeq 0$) means that the square matrix $\mathbf{X}$ is positive definite (semi-definite). $\overline{\tau}(\mathbf{X})$ and $\underline{\tau}(\mathbf{X})$ represent the largest and smallest eigenvalues of $\mathbf{X}$, respectively. $\otimes$ denotes Kronecker product. $\mathbf{I}_n$ is an $n$-dimensional identity matrix and $\mathbf{O}_{n \times m}$ indicates an $(n \times m)$-dimensional zero matrix. $\mathbf{1}_n$ and $\mathbf{0}_n$ denote the $n$-dimensional column vectors with all elements being 1 and 0, respectively.

\subsection{Graph Theory}\label{d2}
Define an undirected graph ${{G}}= \{{{V}},{{E}}\}$ for a multi-agent network, where ${{V}} = \{1,2,...,N\}$ is the set of vertices and ${{E}} =\{e_1,e_2,...,e_{|E|} \} \subseteq \{ (i,j)| i,j \in {V} \hbox{ and } i \neq j\}$ (no self-loop) is the set of edges with $(i,j) \in {E}$ unordered.
Given the index of vertices, the index of edges is determined as follows. For any two distinct edges $e_k = (k_1,k_2) \in E$ and $e_v=(v_1,v_2) \in E$, if $\min\{k_1,k_2\} > \min\{v_1,v_2\}$, then $k > v$, and vice versa. For the case where $\min\{k_1,k_2\} = \min\{v_1,v_2\}$, if $\max\{k_1,k_2\} > \max\{v_1,v_2\}$, then $k > v$, and vice versa.
${G}$ is connected if any two distinct vertices are linked by at least one path. ${{V}}_i =  \{ j | (i,j) \in {E}\}$ represents the neighbor set of agent $i$.
Let ${\mathbf{Q}}$ be the incidence matrix of ${{G}}$ \cite{dimarogonas2010stability}. The $(j,k)$th element of ${\mathbf{Q}}$ is defined by $[\mathbf{Q}]_{jk}=\left\{
\begin{array}{ll}
  1 & \hbox{if }e_k=(j,l)\in E \hbox{ and }j<l \\
  -1 & \hbox{if }e_k=(j,l)\in E\hbox{ and }j>l \\
  0 & \hbox{otherwise}
\end{array}
\right.$.
In addition, define $S_i = \{j| (i,j) \in E, j > i\}$ and ${S}^{\sharp}_i = \{j| (i,j) \in E, j < i\}$. Then, it can be checked that $V_i=S_i \cup S_i^{\sharp}$, $i \in V$.
\subsection{Proximal Mapping}
A proximal mapping of a proper, convex, and closed function $\psi: \mathbb{R}^n \rightarrow (-\infty,+\infty]$ is defined by $\mathrm{prox}^{\alpha}_{\psi} [\mathbf{v}] =  \arg \min_{\mathbf{u}} ( \psi(\mathbf{u}) + \frac{1}{2{\alpha}} \| \mathbf{u} - \mathbf{v}\|^2 )$, ${\alpha}>0$, $\mathbf{v} \in \mathbb{R}^n$ \cite{beck2014fast}.

\subsection{Fenchel Conjugate}\label{de1}
$\psi: \mathbb{R}^n \rightarrow (-\infty,+\infty]$ is a proper function. Then the Fenchel conjugate of $\psi$ is defined by $\psi^{\diamond}(\mathbf{v})= \sup_{\mathbf{u}} \{\mathbf{v}^{\top} \mathbf{u}-\psi(\mathbf{u})\}$, which is convex \cite[Sec. 3.3]{borwein2010convex}.

\begin{Lemma}\label{md}
(Extended Moreau Decomposition \cite[Thm. 6.45]{beck2017first}) $\psi: \mathbb{R}^n \rightarrow (-\infty,+\infty]$ is a proper, convex, and closed function. $\psi^{\diamond}$ is its Fenchel conjugate. Then, for certain $ \mathbf{v} \in \mathbb{R}^n$ and $\alpha >0$, we have $ \mathbf{v} = \alpha \mathrm{prox}^{\alpha^{-1}}_{\psi^{\diamond}} \left[\frac{\mathbf{v}}{\alpha} \right] + \mathrm{prox}^{\alpha}_{\psi} [\mathbf{v}]$.
\end{Lemma}

\begin{Lemma}\label{l1}
\cite[Lemma V.7]{notarnicola2016asynchronous} $\psi: \mathbb{R}^n \rightarrow (-\infty,+\infty]$ is a proper, closed, $\sigma$-strongly convex function and $\psi^{\diamond}$ is its Fenchel conjugate, $\sigma >0$. Then, $\arg \max\limits_{\mathbf{u}} ( \mathbf{v}^{\top} \mathbf{u} - \psi(\mathbf{u})) = \nabla_{\mathbf{v}} \psi^{\diamond}(\mathbf{v})$ and $\nabla_{\mathbf{v}} \psi^{\diamond}(\mathbf{v})$ is Lipschitz continuous with Lipschitz constant $\frac{1}{\sigma}$.
\end{Lemma}

\section{Problem Formulation}\label{sa3}

The problem formulation and relevant assumptions are provided as follows.

Let $F(\mathbf{x})= \sum_{i \in {V}} F_i(\mathbf{x}_i)$ be the global cost function of a multi-agent network ${G}= \{{V},{E}\}$, $\mathbf{x}_i \in \mathbb{R}^M$, $\mathbf{x}=[\mathbf{x}^{\top}_1,...,\mathbf{x}^{\top}_{N}]^{\top} \in \mathbb{R}^{{N}M}$. Agent $i$ maintains a private cost function $F_i(\mathbf{x}_i)= f_i(\mathbf{x}_i)+g_i(\mathbf{x}_i)$. Let $X_i \subseteq \mathbb{R}^M$ be the feasible region of $\mathbf{x}_i$. Then the feasible region of $\mathbf{x}$ can be defined by $X=  X_1 \times X_2 \times ... \times X_{N} \subseteq \mathbb{R}^{NM}$. An affine-constrained optimization problem of $V$ can be given by
\begin{align}
\mathrm{(P1)} \quad   \min\limits_{\mathbf{x} \in X} \quad \sum_{i \in {V}} F_i(\mathbf{x}_i) \quad \hbox{s.t.} \quad \mathbf{A}\mathbf{x} = \mathbf{b}, \nonumber
\end{align}
which is equivalent to
\begin{align}
  \mathrm{(P2)} \quad  \min\limits_{\mathbf{x}} \quad  & \sum_{i \in {V}} (F_i(\mathbf{x}_i) + \mathbb{I}_{X_{i}} (\mathbf{x}_i))
  \quad \hbox{s.t.} \quad  \mathbf{A}\mathbf{x} = \mathbf{b}, \nonumber
\end{align}
with $\mathbf{A} \in \mathbb{R}^{B \times {N}M}$, $\mathbf{b} \in \mathbb{R}^{B}$, $\mathbb{I}_{X_{i}} (\mathbf{x}_i)= \left\{\begin{array}{ll}
                    0 & \hbox{if $\mathbf{x}_i \in X_i$} \\
                    +\infty  & \hbox{otherwise}
                  \end{array}
                  \right.$.

\begin{Assumption}\label{a0}
${G}$ is connected and undirected.
\end{Assumption}

\begin{Assumption}\label{a1}
$f_i: \mathbb{R}^M \rightarrow (-\infty,+\infty]$ and $g_i: \mathbb{R}^M \rightarrow (-\infty,+\infty]$ are proper, convex, and closed extended real-valued functions. In addition, $f_i$ is differentiable and $\sigma_i$-strongly convex, $\sigma_i>0$, $i\in {V}$.
\end{Assumption}

The assumptions in Assumption \ref{a1} are often discussed in composite optimization problems \cite{shi2015proximal,schmidt2011convergence,wang2021composite,chang2014multi,florea2020generalized,beck2014fast,notarnicola2016asynchronous}.

\begin{Assumption}\label{a1+1}
$X_i$ is non-empty, convex, and closed, $i\in {V}$; there exists an $\breve{\mathbf{x}} \in \mathbf{relint}X$ such that $\mathbf{A}\breve{\mathbf{x}} = \mathbf{b}$.
\end{Assumption}

By Assumption \ref{a1+1}, $\mathbb{I}_{X_{i}}$ is proper, convex, and closed \cite{boyd2004convex}, which complies with the assumption on $g_i$.


\section{Distributed Optimization Algorithm Development}\label{sa4}

The proposed algorithm for solving the problem of interest is presented in this section.

\subsection{Dual Problem}\label{sec1}
Decoupling the objective function of Problem (P2) gives
\begin{align}
 \mathrm{(P3)} \quad \min\limits_{\mathbf{x},\mathbf{z}} \quad  & \sum_{i \in {V}} (f_i(\mathbf{x}_i) + (g_i+ \mathbb{I}_{X_{i}}) (\mathbf{z}_i)) \nonumber \\
  \hbox{s.t.} \quad  & \mathbf{A}\mathbf{x} = \mathbf{b}, \quad \mathbf{x}_i = \mathbf{z}_i, \quad \forall i \in {V}, \nonumber
\end{align}
where $\mathbf{z}=[\mathbf{z}^{\top}_1,...,\mathbf{z}^{\top}_{N}]^{\top} \in \mathbb{R}^{NM}$ with $\mathbf{z}_i\in \mathbb{R}^{M}$ being a slack variable. The Lagrangian function of Problem (P3) is
\begin{align}\label{2+12}
   L(\mathbf{x},\mathbf{z},  \bm{\eta},\bm{\mu})= & \sum_{i \in {V}} (f_i(\mathbf{x}_i) + (g_i + \mathbb{I}_{{X}_{i}})(\mathbf{z}_i) \nonumber \\
   & + \bm{\mu}_i^{\top} (\mathbf{x}_i - \mathbf{z}_i))+ \bm{\eta}^{\top} (\mathbf{A}\mathbf{x} - \mathbf{b}) \nonumber \\
   = & \sum_{i \in {V}} (f_i(\mathbf{x}_i) +  \mathbf{x}_i^{\top} ( \mathbf{A}^{\top}_i \bm{\eta} + \bm{\mu}_i) \nonumber \\
& + (g_i+{\mathbb{I}}_{{X}_{i}}) (\mathbf{z}_i)- \mathbf{z}_i^{\top} \bm{\mu}_i) - \mathbf{b}^{\top}\bm{\eta}.
\end{align}
Here, $\bm{\mu}_i \in \mathbb{R}^M$ and $\bm{\eta} \in \mathbb{R}^B $ are the Lagrangian multipliers associated with constraints $\mathbf{x}_i = \mathbf{z}_i$ and $\mathbf{A}\mathbf{x} = \mathbf{b}$, respectively, $\bm{\mu}=[\bm{\mu}^{\top}_1,...,\bm{\mu}^{\top}_{N}]^{\top} \in \mathbb{R}^{{N}M}$, and $\mathbf{A}_i \in \mathbb{R}^{B \times M}$ is the $i$th column sub-block of $\mathbf{A}$ with $\mathbf{A}=[\mathbf{A}_1,...,\mathbf{A}_i,...,\mathbf{A}_{N}]$.

Then the dual function can be obtained by minimizing $L({\mathbf{x}},{\mathbf{z}},\bm{\eta},\bm{\mu})$ with $\mathbf{x}$ and $\mathbf{z}$ \cite{boyd2004convex}, which is $D (\bm{\eta}, \bm{\mu}) = \min\limits_{{\mathbf{x}},{\mathbf{z}}} L(\mathbf{x},\mathbf{z},  \bm{\eta},\bm{\mu})  = \min\limits_{{\mathbf{x}},{\mathbf{z}}} \sum_{i \in {V}} ( f_i(\mathbf{x}_i) -  \mathbf{x}_i^{\top} \mathbf{H}_i \bm{\zeta}_i  +  (g_i + {\mathbb{I}}_{{X}_{i}}) (\mathbf{z}_i)- \mathbf{z}_i^{\top} \mathbf{F} \bm{\zeta}_i - \kappa_i \mathbf{E} \bm{\zeta}_i)  = \sum_{i \in {V}} (- f^{\diamond}_i( \mathbf{H}_i \bm{\zeta}_i ) - \kappa_i \mathbf{E}\bm{\zeta}_i  - ( g_i +{\mathbb{I}}_{{X}_{i}} )^{\diamond}(\mathbf{F}\bm{\zeta}_i) )$, where $\mathbf{H}_i= [-\mathbf{A}^{\top}_i, -\mathbf{I}_M] \in \mathbb{R}^{M \times (M+B)}, \bm{\zeta}_i=   [\bm{\eta}^{\top}, \bm{\mu}_i^{\top}]^{\top} \in \mathbb{R}^{M+B}, \mathbf{F}=  [\mathbf{O}_{M \times B}, \mathbf{I}_M ] \in \mathbb{R}^{M \times (M+B)}, \mathbf{E}=  [\mathbf{b}^{\top},\mathbf{0}^{\top}_{M}] \in \mathbb{R}^{1 \times (M+B)}$, $\sum_{i \in {V}} \kappa_i = 1$, and $(g_i + {\mathbb{I}}_{{X}_{i}}  )^{\diamond}$ denotes the Fenchel conjugate of $g_i + {\mathbb{I}}_{{X}_{i}}$. Hence, the dual problem of Problem (P3) can be formulated as
\begin{align}
\mathrm{(P4)} \quad \min\limits_{\bm{\zeta}} \quad & \sum_{i\in {V}} (f^{\diamond}_i( \mathbf{H}_i \bm{\zeta}_i) + \kappa_i \mathbf{E}\bm{\zeta}_i + (g_i +{\mathbb{I}}_{{X}_{i}} )^{\diamond}(\mathbf{F}\bm{\zeta}_i)), \nonumber
\end{align}
where $\bm{\zeta} = [\bm{\zeta}^{\top}_1 ,...,\bm{\zeta}^{\top}_{N}]^{\top} \in \mathbb{R}^{NB+NM}$. 

We aim to solve Problem (P4) in a distributed manner. In Problem (P4), the variables of $f^{\diamond}_i(\mathbf{H}_i\bm{\zeta}_i)$ are coupled in terms of the common component $\bm{\eta}$ in $\bm{\zeta}_i$, but those of $(g_i +{\mathbb{I}}_{{X}_{i}} )^{\diamond}(\mathbf{F}\bm{\zeta}_i)$ are decoupled since $\mathbf{F}\bm{\zeta}_i=\bm{\mu}_i$. In the following, we define $\bm{\lambda}_i=  [\bm{\theta}_i^{\top}, \bm{\mu}_i^{\top}]^{\top}$ and $\bm{\lambda}= [\bm{\lambda}_1^{\top}, ...,\bm{\lambda}_{N}^{\top}]^{\top}$, where $\bm{\theta}_i  \in \mathbb{R}^{B}$ is the local estimate of $\bm{\eta}$. Then, based on consensus protocol $\bm{\theta}_i=\bm{\theta}_j$, $\forall (i,j) \in E$, Problem (P4) can be equivalently rewritten as
\begin{align}
\mathrm{(P5)} \quad \min\limits_{\bm{\lambda}} \quad & {\Phi}( {\bm{\lambda}}) \nonumber \\
\hbox{s.t.} \quad  & \mathbf{K}\bm{\lambda}_i = \mathbf{K}\bm{\lambda}_j,  \quad \forall (i,j) \in {E}, \label{lb}
\end{align}
where ${\Phi}( {\bm{\lambda}}) = {P}( {\bm{\lambda}}) + {Q}( {\bm{\lambda}}), P (\bm{\lambda}) =  \sum_{i\in {V}} p_i (\bm{\lambda}_i),
Q (\bm{\lambda}) = \sum_{i\in {V}} q_i (\bm{\lambda}_i),
p_i (\bm{\lambda}_i) = f^{\diamond}_i( \mathbf{H}_i \bm{\lambda}_i )  + \kappa_i \mathbf{E}\bm{\lambda}_i, q_i (\bm{\lambda}_i)=  (g_i +{\mathbb{I}}_{{X}_{i}} )^{\diamond}(\mathbf{F}\bm{\lambda}_i), \mathbf{K}=  [\mathbf{I}_B, \mathbf{O}_{B \times M}]$. Constraint (\ref{lb}) ensures the partial consistency among $\bm{\lambda}_i$ in terms of the component $\bm{\theta}_i$, i.e., $\bm{\theta}_i = \mathbf{K}\bm{\lambda}_i$.
In addition, (\ref{lb}) can be represented by $\mathbf{K}\bm{\lambda}_i = \mathbf{K}\bm{\lambda}_j$, $\forall i \in V, j\in S_i$ (if $S_i \neq \emptyset$), which can be written into a compact form $\mathbf{M}\bm{\lambda} = \mathbf{0}$, where ${\mathbf{M}} = \mathbf{Q}^{\top} \otimes \mathbf{K}$ \cite{dimarogonas2010stability}.
Then, Problem (P5) is equivalent to
\begin{align}
\mathrm{(P6)} \quad \min\limits_{ {\bm{\lambda}} } \quad & {\Phi}( {\bm{\lambda}})  \quad
\hbox{s.t.} \quad  \mathbf{M}\bm{\lambda} = \mathbf{0}. \nonumber
\end{align}
Let $\bm{\lambda}^*= [\bm{\lambda}_1^{*\top}, ...,\bm{\lambda}_{N}^{*\top}]^{\top}$ be the optimal solution to Problem (P6) with $\bm{\lambda}^*_i=[\bm{\theta}^{*\top}_i,\bm{\mu}^{*\top}_i]^{\top}$.

The Lagrangian function of Problem (P6) can be given by
 \begin{align}\label{fp1}
& \mathcal{L} (\bm{\lambda},\bm{\xi}) =  P(\bm{\lambda})+ Q(\bm{\lambda}) +  \bm{\xi}^{\top} {\mathbf{M}} \bm{\lambda},
\end{align}
where $\bm{\xi} =[\bm{\xi}^{\top} _1,...,\bm{\xi}^{\top}_{N}]^{\top}$, $\bm{\xi}_{i} = [ \bm{\xi}^{\top}_{ij_1},...,\bm{\xi}^{\top}_{ij_{|S_{i}|}}]^{\top}$ with $\bm{\xi}_{ij_l} \in \mathbb{R}^B$ being the Lagrangian multiplier associated with the constraint $\mathbf{K}\bm{\lambda}_i - \mathbf{K}\bm{\lambda}_{j_l} = \mathbf{0}$, $j_l \in S_i$, and $l \in \{1,2,...,|S_i|\}$ is the index. Let ${C}$ be the set of the saddle points of $\mathcal{L} (\bm{\lambda},\bm{\xi})$. Then any saddle point $(\bm{\lambda}^*,\bm{\xi}^*) \in {C}$ satisfies \cite{rockafellar1970convex}
\begin{align}\label{sad}
 \mathcal{L} (\bm{\lambda},\bm{\xi}^*) \geq  \mathcal{L} (\bm{\lambda}^*,\bm{\xi}^*) \geq  \mathcal{L} (\bm{\lambda}^*,\bm{\xi}).
\end{align}
We aim to seek a saddle point of $\mathcal{L}(\bm{\lambda},\bm{\xi}) $, which
can be characterized by Karush-Kuhn-Tucker (KKT) conditions \cite{boyd2004convex}
\begin{align}
 & \mathbf{0}\in \nabla_{\bm{\lambda}} P(\bm{\lambda}^*) + \partial_{\bm{\lambda}} Q(\bm{\lambda}^*)
 + \mathbf{M}^{\top} \bm{\xi}^*,    \label{k1} \\
&  \mathbf{M} \bm{\lambda}^* = \mathbf{0}. \label{k2}
\end{align}

\subsection{DDPG Algorithm}

Based on the previous discussion, the DDPG algorithm for solving Problem (P6) is designed as
\begin{align}
 \bm{\lambda}^{t+1} = & \mathrm{prox}^{c}_{Q}\big[\bm{\lambda}^t - c  ( \nabla_{\bm{\lambda}} P (\bm{\lambda}^t)   + {\mathbf{M}}^{\top} \bm{\xi}^t + \gamma \mathbf{M}^{\top} \mathbf{M} \bm{\lambda}^t) \big], \label{f1} \\
 \bm{\xi}^{t+1}= & \bm{\xi}^t + \gamma {\mathbf{M}} \bm{\lambda}^{t+1}, \label{f2}
\end{align}
which means
\begin{align}
 \bm{\lambda}_i^{t+1} = & \mathrm{prox}^{c}_{q_i} \bigg[\bm{\lambda}_i^t - c  \bigg( \nabla_{\bm{\lambda}_i} p_i (\bm{\lambda}_i^t)  + \sum_{j \in {S}_i} \mathbf{K}^{\top}\bm{\xi}_{ij}^t
  \nonumber \\
  - \sum_{j \in {S}^{\sharp}_i} & \mathbf{K}^{\top}  \bm{\xi}_{ji}^t  + \gamma \sum_{j \in {V}_i} \mathbf{K}^{\top} \mathbf{K} (\bm{\lambda}_i^t - \bm{\lambda}_j^t)\bigg) \bigg] , \forall i \in V, \label{f3}   \\
 \bm{\xi}_{ij}^{t+1}  = & \bm{\xi}_{ij}^t  + \gamma \mathbf{K} (\bm{\lambda}_{i}^{t+1} - \bm{\lambda}_j^{t+1}), \quad \forall i \in V, j \in S_i,
 \label{f4}
\end{align}
due to the separability of $P$ and $Q$, $t \in \mathbb{N}$. $c,\gamma>0$ are step sizes.

\begin{Remark}
As seen in (\ref{f1}), we employ $\gamma \mathbf{M}^{\top} \mathbf{M} \bm{\lambda}^t$ in the updating law to improve the stability of the algorithm. Then, as shown in (\ref{f3}) and (\ref{f4}), each agent only communicates with its neighbors (i.e., a distributed computation manner). To this end, the incidence matrix $\mathbf{Q}$ plays a key role in formulating the desired structure of $\gamma \mathbf{M}^{\top} \mathbf{M} \bm{\lambda}^t$.  As another feature, the update of $\bm{\xi}_{i}$ requires agent $i$ to communicate only with the agents in set $S_i$ rather than all its neighbors in $V_i$, which can reduce the communication burden in each updating round.
\end{Remark}

\begin{Remark}\label{rm2}
If the structure of $f_i$ is complicated and $\nabla p_i$ cannot be obtained efficiently, we can implement (\ref{f3}) by computing
\begin{align}\label{r7+2}
\nabla_{\bm{\lambda}_i} p_i & (\bm{\lambda}_i^t) = \nabla_{\bm{\lambda}_i} f^{\diamond}_i( \mathbf{H}_i \bm{\lambda}_i ^t) + \kappa_i \mathbf{E}^{\top} \nonumber \\
= & \mathbf{H}_i^{\top} \nabla_{\mathbf{H}_i\bm{\lambda}_i}f_i^{\diamond} ( \mathbf{H}_i \bm{\lambda}_i ^t) + \kappa_i \mathbf{E}^{\top} \nonumber \\
= & \mathbf{H}_i^{\top} \arg \max_{\mathbf{u}} ((\mathbf{H}_i\bm{\lambda}_i ^t)^{\top} \mathbf{u} -f_i(\mathbf{u})) + \kappa_i \mathbf{E}^{\top}
\end{align}
based on Lemma \ref{l1}. In this case, the assumption on the smooth parts can be extended to non-smooth cases but can be with a higher computational complexity by computing (\ref{r7+2}).
\end{Remark}

In the following, we will discussed how to recover the optimal primal solution. Note that the optimal solution $\bm{\lambda}^*$ to Problem (P6) is also the optimal solution to Problem (P4) due to the equivalence. Then the optimal primal solution to Problem (P3) can be characterized by the saddle point property
\begin{align}\label{lb-1}
L (\mathbf{x}^*, \mathbf{z}^*, \bm{\eta},\bm{\mu}) \leq {L (\mathbf{x}^*, \mathbf{z}^*, \bm{\eta}^*,\bm{\mu}^*) \leq L (\mathbf{x}, \mathbf{z}, \bm{\eta}^*,\bm{\mu}^*)},
\end{align}
where $\bm{\eta}^*=\bm{\theta}_{i}^*$ since $\bm{\theta}_{i}^*$ is the consensual local estimate of $\bm{\eta}^*$, $i \in V$. Therefore, the optimal primal solution can be obtained by the second inequality in (\ref{lb-1}), which gives
\begin{align}\label{ala}
\mathbf{x}^* = \arg \min_{\mathbf{x}} L (\mathbf{x}, \mathbf{z}, \bm{\eta}^*,\bm{\mu}^*).
\end{align}
Note that $\mathbf{x}^*$ is unique since $f_i$ is strongly convex. Then by omitting some constant terms in $L$ and decomposing $\mathbf{x}^*$, (\ref{ala}) can be rewritten as $
\mathbf{x}_i^* = \arg \min_{\mathbf{x}_i}f_i(\mathbf{x}_i) +  \mathbf{x}_i^{\top} ( \mathbf{A}^{\top}_i \bm{\eta}^* + \bm{\mu}^*_i)  = \arg \min_{\mathbf{x}_i}f_i(\mathbf{x}_i) +  \mathbf{x}_i^{\top} ( \mathbf{A}^{\top}_i \bm{\theta}_i^* + \bm{\mu}^*_i) = \arg \min_{\mathbf{x}_i}f_i(\mathbf{x}_i) - \mathbf{x}_i^{\top} \mathbf{H}_i \bm{\lambda}^*_i$, which can be completed by agent $i$ locally.

The detailed computation procedure of the DDPG algorithm is stated in Algorithm \ref{ax}.

\begin{algorithm}
\caption{DDPG Algorithm}\label{ax}
\begin{algorithmic}[1]
\State Initialize $\bm{\lambda}^0$, $\bm{\xi}^0$. Determine step sizes $c,\gamma>0$.
\For {$t= 0,1,2,...$}
\For {$i= 1,2,...,{N}$} (in parallel)
\State Update $\bm{\lambda}_i^{t+1}$ based on  (\ref{f3}).
\State Update $\bm{\xi}_{ij}^{t+1}$ based on  (\ref{f4}), $\forall j \in S_i$.
\EndFor
\EndFor
\State Obtain the outputs $\bm{\lambda}_i^{\mathrm{out}}$ and $\bm{\xi}_{ij}^{\mathrm{out}}$ under certain convergence criterion and calculate the primal solution by $\mathbf{x}_i^{\mathrm{out}} = \arg \min_{\mathbf{x}_i}f_i(\mathbf{x}_i) - \mathbf{x}_i^{\top} \mathbf{H}_i \bm{\lambda}^{\mathrm{out}}_i$, $\forall i \in V, j\in S_i$.
\end{algorithmic}
\end{algorithm}

\section{Main Result}\label{sa5}

The convergence analysis and computational complexity analysis of the DDPG algorithm are conducted in this section.


\subsection{Convergence Analysis}

\begin{Lemma}\label{lam2}
Based on Assumption \ref{a1}, the Lipschitz constant of $\nabla_{\bm{\lambda}_i} p_i (\bm{\lambda}_i)$ is given by $h_i= \frac{ \| \mathbf{H}_i \|^2 }{\sigma_i}$, $i \in V$.
\end{Lemma}

See the proof in Appendix \ref{lam2p}.


\begin{Theorem}\label{la1}
Suppose that Assumptions \ref{a0}-\ref{a1+1} hold. If the step sizes satisfy
\begin{align}\label{37}
c^{-1} \geq & {h} +  \gamma \overline{\tau} \left( \mathbf{M}^{\top} \mathbf{M} \right)
\end{align}
with ${h}=\max \{h_i\}_{i\in V} $, then the states generated by Algorithm \ref{ax} converge to a primal-dual solution to Problem (P6).
\end{Theorem}

See the proof in Appendix \ref{la1p}.

\begin{Theorem}\label{th1}
Suppose that Assumptions \ref{a0}-\ref{a1+1} hold and the step sizes are selected based on (\ref{37}). By Algorithm \ref{ax}, for certain $(\bm{\lambda}^*,\bm{\xi}^*) \in {C}$, we have $ | \Phi (\overline{\bm{\lambda}}^{T+1}) - \Phi(\bm{\lambda}^*) | \leq \mathcal{O}\left(\frac{1}{T}\right)$ and $\| \bm{\xi}^*  \| \| \mathbf{M} \overline{\bm{\lambda}}^{T+1} \| \leq  \mathcal{O}\left(\frac{1}{T}\right)$, where $\overline{\bm{\lambda}}^{T+1} = \frac{1}{T+1}\sum_{t=0}^T \bm{\lambda}^{t+1}$ and $T \in \mathbb{N}_+$.
\end{Theorem}

See the proof in Appendix \ref{th1p}.

\subsection{Computational Complexity Analysis with Simple-Structured Cost Functions}\label{c3}

In this work, the concept of computational complexity measures the total amount of basic operations required, which is dominated by the iteration complexity and computational cost per iteration \cite{bubeck2015convex}.

In the following, we first discuss the iteration complexity of (\ref{f3}). Specifically, to apply (\ref{f3}), one needs to compute (i) $\nabla p_i$ and (ii) the proximal mapping of $q_i$, $i\in {V}$. For (i), $\nabla p_i$ can be efficiently obtained given that $f_i$ is simple-structured and, consequently, $\nabla f^{\diamond}_i$ can be analytically derived, e.g., $f_i$ is a quadratic function \cite[Sec. 3.3.1]{boyd2004convex}. For (ii), some feasible methods for different cases are introduced as follows.

{\em{1) Case 1:}} If the proximal mapping of $g_i+\mathbb{I}_{{X}_{i}}$ can be easily obtained,\footnote[3]{This case is based on that $g_i+\mathbb{I}_{{X}_{i}}$ is with certain simple structure, which is often the assumption in the works on proximal gradient method. See some frequently used formulas in \cite[Sec. 6.3]{beck2017first} and applications in \cite[Sec. 7]{parikh2014proximal}.} we have $\mathrm{prox}^{c}_{q_i} =  \mathrm{prox}^{c}_{(g_i+{\mathbb{I}}_{{X}_{i}})^{\diamond}}$, where $\mathrm{prox}^{c}_{(g_i+{\mathbb{I}}_{{X}_{i}})^{\diamond}}$ can be obtained by calculating $\mathrm{prox}^{ c^{-1}}_{g_i+{\mathbb{I}}_{{X}_{i}}}$ based on Lemma \ref{md}. Then, by decomposing $\bm{\lambda}_i^{t+1}$, (\ref{f3}) can be modified into $\bm{\theta}_i^{t+1} = \bm{\theta}_i^t - c ( \nabla_{\bm{\theta}_i} p_i (\bm{\lambda}_i^t) + \sum_{j \in {S}_i} \bm{\xi}_{ij}^t  - \sum_{j \in {S}^{\sharp}_i} \bm{\xi}_{ji}^t + \gamma \sum_{j \in {V}_i} (\bm{\theta}_i^t - \bm{\theta}_j^t) )$ and
\begin{align}
& \bm{\mu}_i^{t+1} = \mathrm{prox}^{c}_{{q}_i} \big[ \bm{\mu}_i^t - c \nabla_{\bm{\mu}_i} p_i (\bm{\lambda}_i^t) \big] \nonumber \\
& =  \bm{\mu}_i^t  - c \nabla_{\bm{\mu}_i} p_i (\bm{\lambda}_i^t) - c \mathrm{prox}^{c^{-1}}_{{q}^{\diamond}_{i}} \left[\frac{\bm{\mu}_i^t - c \nabla_{\bm{\mu}_i} p_i (\bm{\lambda}_i^t)}{c} \right], \label{r7+1}
\end{align}
where $q^{\diamond}_{i}=(g_i +{\mathbb{I}}_{{X}_{i}} )^{\diamond \diamond}=g_i +{\mathbb{I}}_{{X}_{i}} $ due to the convexity and lower semi-continuity of $g_i +{\mathbb{I}}_{{X}_{i}}$ and $(g_i +{\mathbb{I}}_{{X}_{i}} )^{\diamond \diamond}$ is the biconjugate of $g_i +{\mathbb{I}}_{{X}_{i}}$ \cite[Sec. 3.3.2]{boyd2004convex}. Then, the calculation of the proximal mapping of $(g_i +{\mathbb{I}}_{{X}_{i}})^{\diamond} $ can be avoided as shown in (\ref{r7+1}), which can reduce the iteration complexity if the proximal mapping of $g_i +{\mathbb{I}}_{{X}_{i}}$ is easier to obtain. For instance, in an LASSO problem with penalty $g_i(\mathbf{x}_i)= \| \mathbf{x}_i \|_1$ and $X_i = \mathbb{R}^M$, the proximal mapping of $l_1$-norm is a soft thresholding operator with analytical solution \cite[Sec. 6.3]{beck2017first}. In addition, if $g_i=0$, (\ref{r7+1}) can be written as $\bm{\mu}_i^{t+1} = \bm{\mu}_i^t - c \nabla_{\bm{\mu}_i} p_i (\bm{\lambda}_i^t) -  c \mathrm{prox}^{c^{-1}}_{\mathbb{I}_{{X}_{i}}} \left[\frac{\bm{\mu}_i^t - c \nabla_{\bm{\mu}_i} p_i (\bm{\lambda}_i^t)}{c}\right]
 = \bm{\mu}_i^t - c \nabla_{\bm{\mu}_i} p_i (\bm{\lambda}_i^t) -  c \Pi_{{X_i}} \left[\frac{\bm{\mu}_i^t - c \nabla_{\bm{\mu}_i} p_i (\bm{\lambda}_i^t)}{c}\right]$, where $\Pi_{{X_i}}[\cdot]$ is a Euclidean projection onto $X_i$ \cite[Sec. 1.2]{parikh2014proximal}.

{\em{2) Case 2:}} Take the advantage of the structure of $g_i$ in some specific problems. For example, consider a regularization problem, where the penalty is a Euclidean $e$-norm: $g_i(\mathbf{x}_i)= \| \mathbf{x}_i \|_e$, $X_i = \mathbb{R}^M$. Then we can have $q_i (\bm{\lambda}_i) = g^{\diamond}_i( \bm{\mu}_i)
= \mathbb{I}_{W_i} (\bm{\mu}_i) = \left\{\begin{array}{ll}
                    0 & \hbox{if $\bm{\mu}_i \in W_i$} \\
                    +\infty  & \hbox{otherwise}
                  \end{array}
                  \right. = \left\{\begin{array}{ll}
                    0 & \hbox{if $\bm{\lambda}_i \in Y_i$ } \\
                    +\infty  & \hbox{otherwise}
                  \end{array}
                  \right. = \mathbb{I}_{Y_i} (\bm{\lambda}_i)$,
where $W_i=\{ \mathbf{v}  \in \mathbb{R}^{M} | \| \mathbf{v} \|^*_{e} \leq 1 \}$ (convex zone) with $\| \cdot \|^*_{e}$ being the dual norm of $\| \cdot \|_e$, and $Y_i =  \mathbb{R}^B \times W_i$. The second equality holds by computing the conjugate of a norm \cite[Sec. 3.3.1]{boyd2004convex}. Then, in (\ref{f3}), the proximal mapping of $q_i$ is a Euclidean projection onto $Y_i$ \cite[Sec. 1.2]{parikh2014proximal}.

{\em{3) Case 3:}} If $q_i$ is with certain complicated structure, as a general method, we can construct a strongly convex non-smooth $g_i$ (e.g., shift a strongly convex component of $f_i$ to $g_i$). Then, rewrite (\ref{f3}) by the definition of proximal mapping, which gives $\bm{\lambda}_i  ^{t+1} =  \arg \min_{\bm{\lambda}_i} (q_i(\bm{\lambda}_i) + \frac{1}{2c} \| \bm{\lambda}_i - \bm{\lambda}_i^t  + c  ( \nabla_{\bm{\lambda}_i} p_i (\bm{\lambda}_i^t) + \sum_{j \in {S}_i} \mathbf{K}^{\top}\bm{\xi}_{ij}^t -  \sum_{j \in {S}^{\sharp}_i}  \mathbf{K}^{\top}  \bm{\xi}_{ji}^t + \gamma \sum_{j \in {V}_i}\mathbf{K}^{\top}  \mathbf{K} (\bm{\lambda}_i^t - \bm{\lambda}_j^t) ) \|^2 )$.
To obtain the result, one can utilize a gradient descent method by computing $\nabla_{\bm{\lambda}_i} q_i(\bm{\lambda}_i) = \nabla_{\bm{\lambda}_i}(g_i +{\mathbb{I}}_{{X}_{i}})^{\diamond}(\mathbf{F}\bm{\lambda}_i) = \mathbf{F}^{\top} \nabla_{\mathbf{F}\bm{\lambda}_i}(g_i +{\mathbb{I}}_{{X}_{i}})^{\diamond}(\mathbf{F}\bm{\lambda}_i) = \mathbf{F}^{\top} \arg \max_{\mathbf{u}} ((\mathbf{F}\bm{\lambda}_i)^{\top} \mathbf{u} -(g_i + \mathbb{I}_{X_i})(\mathbf{u}))$ with the help of Lemma \ref{l1}. In this case, the update of $\bm{\lambda}_i$ requires an inner-loop optimization to compute the gradient of $q_i$, which can be completed by agent $i$ locally.

In Cases 1 and 2, (\ref{f3}) can only involve some simple operations (e.g., addition, multiplication, and Euclidean projection with iteration complexity $\mathcal{O}(1)$) without any costly inner-loop optimization of primal or other auxiliary variables, which results in an overall iteration complexity $\mathcal{O}\left(\frac{1}{\epsilon}\right)$ of (\ref{f3}) with an ergodic convergence error $\epsilon$ in dual function value (see Thm. \ref{th1}). In addition, note that the computational cost per iteration is linear in the dimension of $\bm{\lambda}$ and $\bm{\xi}$ \cite{bubeck2015convex}, then the overall computational complexity of the DDPG algorithm can be $ \mathcal{O}\left(\frac{NM+NB}{\epsilon}\right)+\mathcal{O}\left(\frac{B\sum_{i\in V}|S_i|}{\epsilon}\right) $.


\section{Numerical Result}\label{sa6}

\begin{figure}
  \centering
  \includegraphics[height=2.5cm,width=4.2cm]{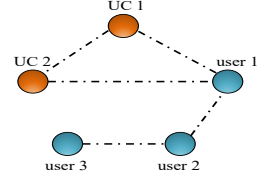}\\
  \caption{Communication typology of the market.}\label{ty}
\end{figure}

In this section, we demonstrate the feasibility of the DDPG algorithm by solving a social welfare optimization problem in an electricity market, where the overall cost function is the sum of the objective function of all the agents \cite{wang2022social}.

\begin{figure*}[t]
\centering
\subfigure[Values of $\bm{\theta}$ and $\bm{\mu}$.]{
\begin{minipage}[t]{0.33\linewidth}
\centering
\includegraphics[width=6.5cm,height=3.5cm]{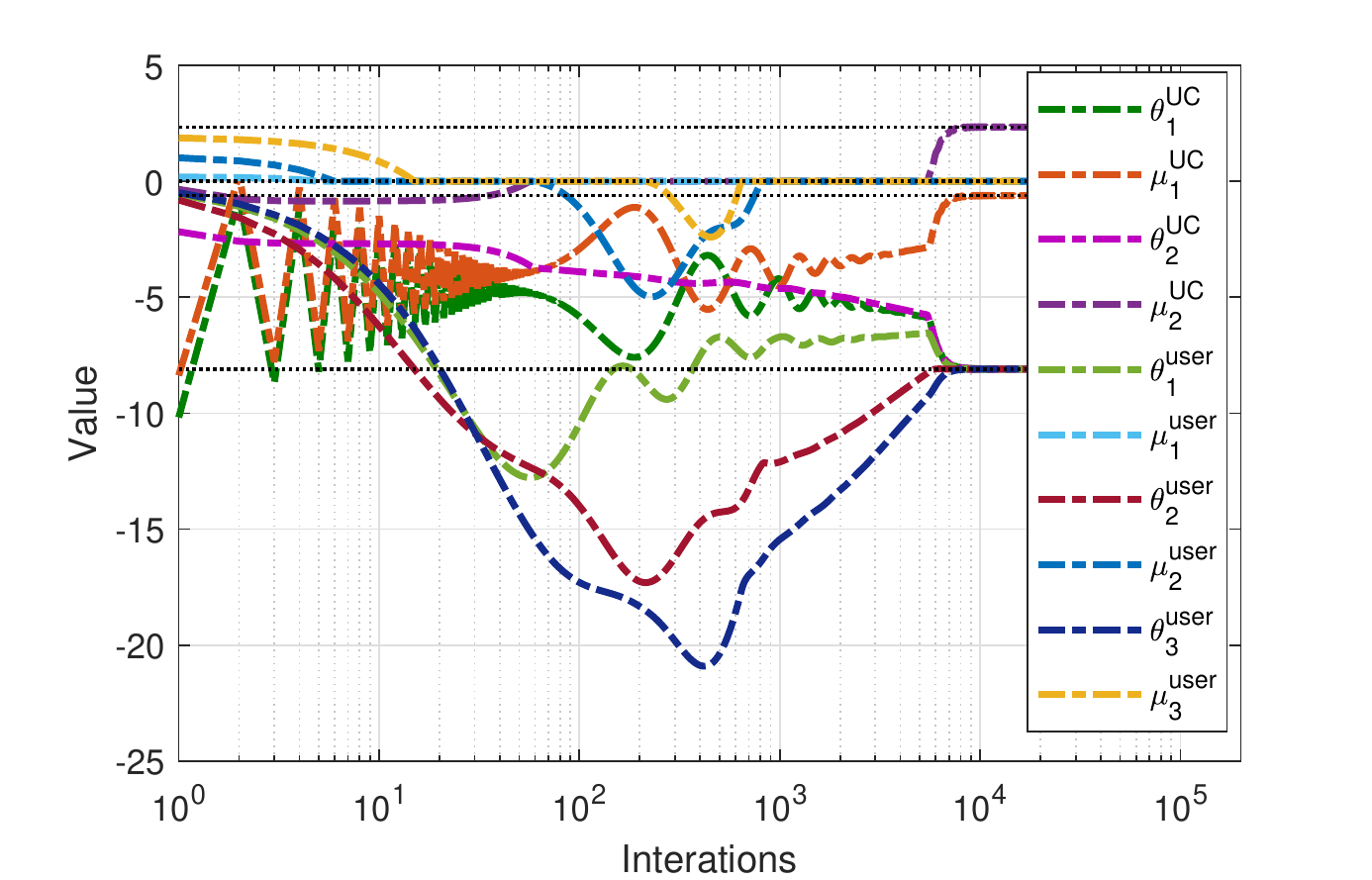}
\end{minipage}%
}%
\subfigure[Values of $\bm{\xi}$.]{
\begin{minipage}[t]{0.33\linewidth}
\centering
\includegraphics[width=6.5cm,height=3.5cm]{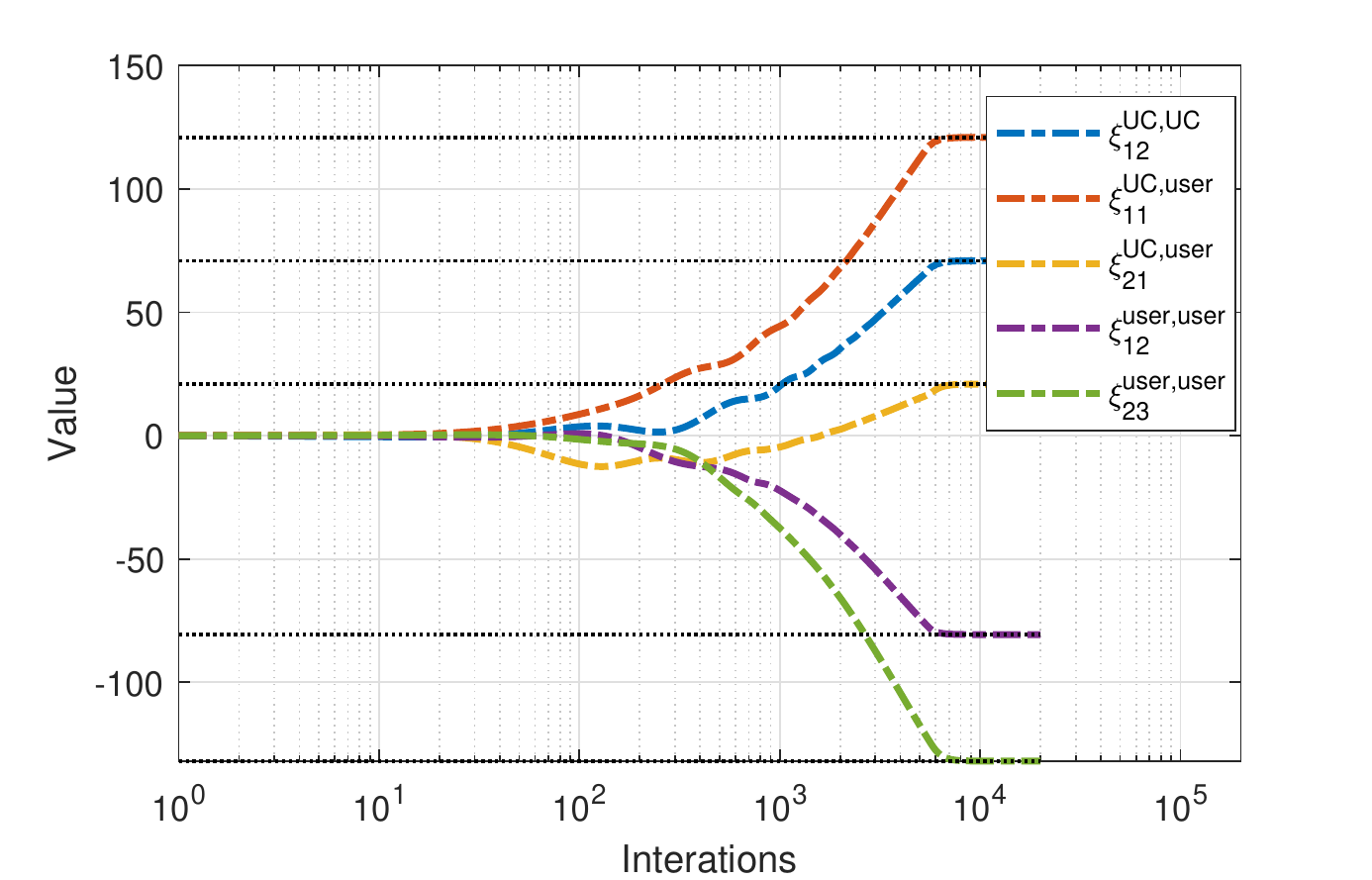}
\end{minipage}%
}%
\subfigure[Values of ${\Phi}( {\bm{\lambda}})$.]{
\begin{minipage}[t]{0.33\linewidth}
\centering
\includegraphics[width=6.5cm,height=3.5cm]{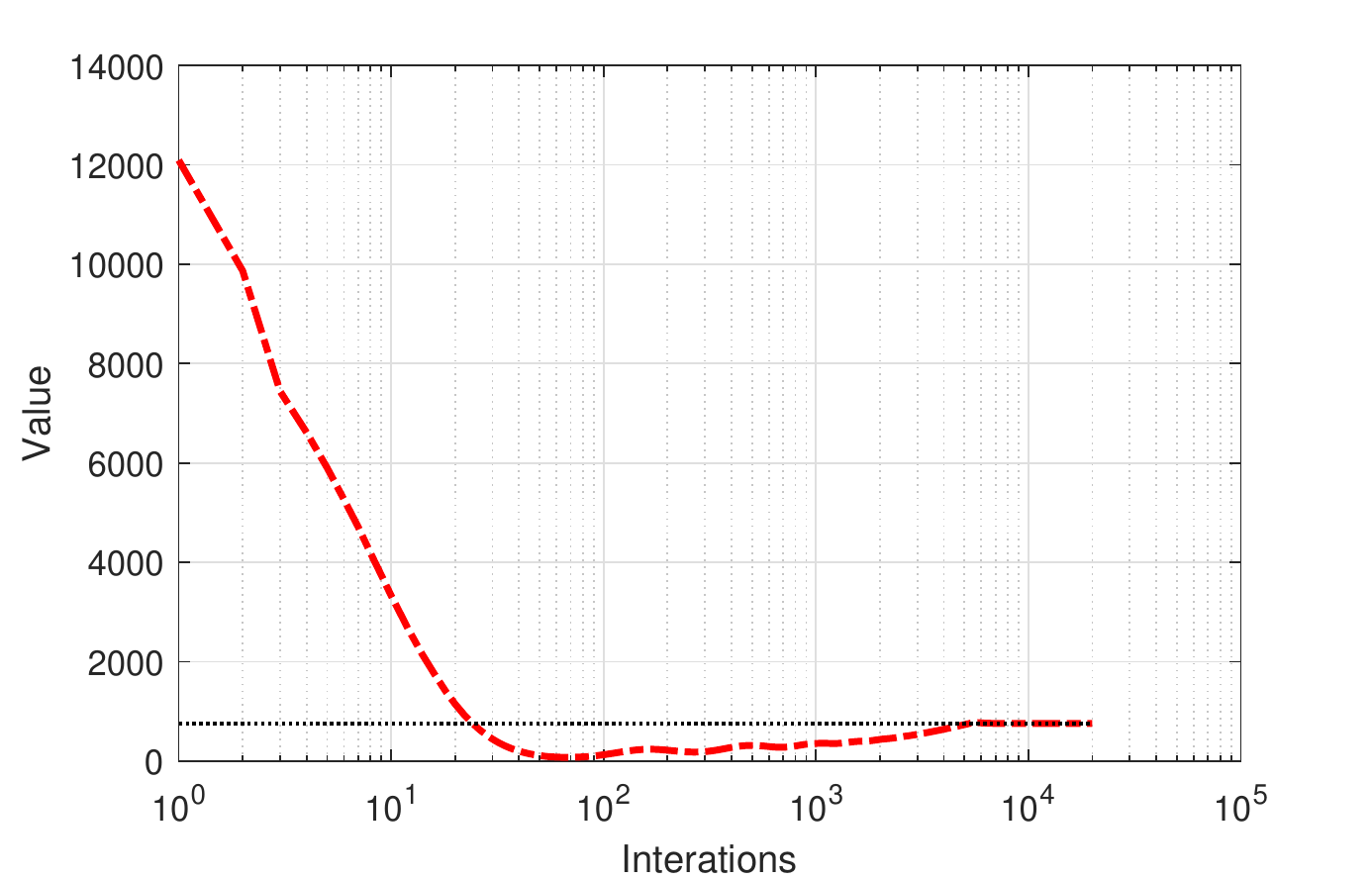}
\end{minipage}%
}%
\caption{Simulation result.}\label{g2}
\end{figure*}

\subsection{Simulation Setup}

Define ${V}_{\mathrm{UC}}$ and ${V}_{\mathrm{user}}$ as the sets of utility companies (UCs) and users, respectively. Define $\mathbf{x} = [x^{\mathrm{UC}}_{1},..., x^{\mathrm{UC}}_{|{V}_{\mathrm{UC}}|}, x^{\mathrm{user}}_{1},...,x^{\mathrm{user}}_{|{V}_{\mathrm{user}}|}]^{\top}$, where $x^{\mathrm{UC}}_{i}$ is the energy generation quantity of UC $i$ and $x^{\mathrm{user}}_{j}$ is the demand of user $j$, $i \in {V}_{\mathrm{UC}}$, $j \in {V}_{\mathrm{user}}$. $\phi_i(x^{\mathrm{UC}}_{i})$ and $\omega_j(x^{\mathrm{user}}_{j})$ are the cost function of UC $i$ and the utility function of user $j$, respectively. Then the social welfare optimization problem of the market can be formulated as
\begin{align}
\mathrm{(P7)} \quad \min_{\mathbf{x}} \quad & \sum_{i \in {V}_{\mathrm{UC}}} \phi_i(x^{\mathrm{UC}}_{i}) - \sum_{j \in {V}_{\mathrm{user}}} \omega_j(x^{\mathrm{user}}_{j}) \nonumber\\
  \hbox{s.t.} \quad & \sum_{i \in {V}_{\mathrm{UC}}} x^{\mathrm{UC}}_{i} = \sum_{j \in {V}_{\mathrm{user}}} x^{\mathrm{user}}_{j}, \label{s1} \\
  & x^{\mathrm{UC}}_{i} \in X^{\textrm{UC}}_i, x^{\mathrm{user}}_{j} \in X^{\textrm{user}}_j, \forall i \in  {V}_{\mathrm{UC}}, j \in  {V}_{\mathrm{user}}, \nonumber
\end{align}
where $\phi_i(x^{\mathrm{UC}}_{i}) = \delta_i (x^{\mathrm{UC}}_{i})^2 + \varsigma_i x^{\mathrm{UC}}_{i} +\beta_i$, $\omega_j(x^{\mathrm{user}}_{j}) =
\left\{\begin{array}{ll}
\chi_j x^{\mathrm{user}}_{j} - \pi_j (x^{\mathrm{user}}_{j})^2 & x^{\mathrm{user}}_{j} \leq \frac{\chi_j}{2\pi_j} \\
  \frac{\chi_j^2}{4 \pi_j} & x^{\mathrm{user}}_{j} > \frac{\chi_j}{2\pi_j}
\end{array}
\right.$, with $\delta_i,\varsigma_i,\beta_i,\chi_j$, and $\pi_j$ being parameters, whose values are set in Table I \cite{pourbabak2017novel}. (\ref{s1}) is the supply-demand balance constraint. $X^{\textrm{UC}}_i = [0,x^{\mathrm{UC}}_{i,\mathrm{max}}]$ and $X^{\textrm{user}}_j = [0,x^{\mathrm{user}}_{j,\mathrm{max}}]$ are local constraints with $x^{\mathrm{UC}}_{i,\mathrm{max}},x^{\mathrm{user}}_{j,\mathrm{max}}>0$. Define $ {\mathbf{A}} =  [\mathbf{1}^{\top}_{|{V}_{\mathrm{UC}}|}, -\mathbf{1}^{\top}_{|{V}_{\mathrm{user}}|}]$. Then (\ref{s1}) is equivalent to ${\mathbf{A}}\mathbf{x} = 0$.

\begin{table}\label{tm2}
\caption{Parameters of UCs and energy users}
\label{tab2}
\begin{center}
\begin{tabular}{p{4mm}p{7mm}p{6mm}p{4mm}p{7mm}p{6mm}p{8mm}p{7mm}}
\bottomrule
& \multicolumn{3}{l}{$ \quad \quad \quad\quad$ UCs} & \multicolumn{3}{l}{$\quad \quad \quad \quad \quad $ Users} \\
\hline
$i/j$ & $\delta_i$ & $\varsigma_i$ & $\beta_i$ & $x^{\mathrm{UC}}_{i,\mathrm{max}}$ & $\chi_j$ &  $\pi_j$  & $x^{\mathrm{user}}_{j,\mathrm{max}}$ \\
\hline
1 & 0.0031& 8.71 & 0 &150 & 17.17 & 0.0935 & 91.79 \\
2 & 0.0074 & 3.53 & 0 &150 & 12.28 & 0.0417 & 147.29 \\
3 & - & - & - &  - &18.42 & 0.1007 & 91.41 \\
\bottomrule
\end{tabular}
\end{center}
\end{table}


To show the performance of Algorithm \ref{ax}, we consider the communication typology shown in Fig. \ref{ty}. Similar to the derivation procedure of (\ref{2+12}), the Lagrangian function of Problem (P7) can be obtained as $L(\mathbf{x},\mathbf{z}, {\eta},\bm{\mu})=  \sum_{i \in {V}_{\mathrm{UC}}} (\phi_i(x^{\mathrm{UC}}_{i})+ \mathbb{I}_{X^{\textrm{UC}}_i} (z^{\mathrm{UC}}_{i})) + \sum_{j \in {V}_{\mathrm{user}}} ( -\omega_j(x^{\mathrm{user}}_{j})+ \mathbb{I}_{X^{\textrm{user}}_j} (z^{\mathrm{user}}_{j}))  + {\eta} {\mathbf{A}}\mathbf{x} + \sum_{i \in {V}_{\mathrm{UC}}}{\mu}^{\mathrm{UC}}_i ({x}^{\mathrm{UC}}_i - {z}^{\mathrm{UC}}_i)  + \sum_{j \in {V}_{\mathrm{user}}} {\mu}^{\mathrm{user}}_j ( {x}^{\mathrm{user}}_j - {z}^{\mathrm{user}}_j)$, where $\mathbf{z}=[ z^{\mathrm{UC}}_1,..., z^{\mathrm{UC}}_{\mid {V}_{\mathrm{UC}}\mid}, z^{\mathrm{user}}_1,..., z^{\mathrm{user}}_{\mid {V}_{\mathrm{user}}\mid}]^{\top}$ is a slack variable and both $\eta$ and $\bm{\mu}=[ \mu^{\mathrm{UC}}_1,..., \mu^{\mathrm{UC}}_{\mid {V}_{\mathrm{UC}}\mid}, {\mu}^{\mathrm{user}}_1,...,$ ${\mu}^{\mathrm{user}}_{\mid {V}_{\mathrm{user}}\mid}]^{\top}$ are dual variables. Define $\bm{\theta}= [ {\theta}^{\mathrm{UC}}_1,..., {\theta}^{\mathrm{UC}}_{\mid {V}_{\mathrm{UC}}\mid},$ $ {\theta}^{\mathrm{user}}_1,...,{\theta}^{\mathrm{user}}_{\mid {V}_{\mathrm{user}}\mid}]^{\top}$, which contains the local estimates of ${\eta}$.
Let $\bm{\xi} = [{\xi}^{\mathrm{UC},\mathrm{UC}}_{12}, {\xi}^{\mathrm{UC},\mathrm{user}}_{11}, {\xi}^{\mathrm{UC},\mathrm{user}}_{21}, {\xi}^{\mathrm{user},\mathrm{user}}_{12}, {\xi}^{\mathrm{user},\mathrm{user}}_{23}]^{\top}$ be the Lagrangian multiplier defined in (\ref{fp1}). For instance, ${\xi}^{\mathrm{UC},\mathrm{user}}_{21}$ is the multiplier associated with constraint $\theta^{\mathrm{UC}}_2= \theta^\mathrm{user}_1$. With some direct calculations, the optimal solution to Problem (P7) is given by $\mathbf{x}^*=[0,150,48.5,50.2,51.3]^{\top}$. 

\subsection{Simulation Result}

The simulation result is shown in Figs. \ref{g2}-(a) to \ref{g2}-(c). Fig. \ref{g2}-(a) depicts the trajectory of dual variables $\bm{\theta}$ and $\bm{\mu}$. It can be seen that all the elements in $\bm{\theta}$ converge to $\eta^*=-8.1$ while $\bm{\mu}$ converges to $\bm{\mu}^*= [-0.61,2.34,0,0,0]^{\top}$. One can check that the optimal primal solution at the saddle point of $L$ is $\mathbf{x}^*= \arg \min_{\mathbf{x}} L( \mathbf{x}, \mathbf{z}, \eta^*,\bm{\mu}^*) = [0,150,48.5,50.2,51.3]^{\top}$, which means that the lower bound and upper bound of $x^{\mathrm{UC}}_{1}$ and $x^{\mathrm{UC}}_{2}$ are activated, respectively, while other variables reach interior optimal solutions. Fig. \ref{g2}-(b) depicts the trajectory of $\bm{\xi}$. Fig. \ref{g2}-(c) shows that the value of dual function ${\Phi}( {\bm{\lambda}})$ (as defined in Problem (P6)) converges to approximately 756.53.

\section{Conclusion}\label{sa7}

In this work, we considered solving a composite DOP with affine coupling constraints. A DDPG algorithm was proposed by resorting to the dual problem. Compared with the existing research works with similar problem setups, we showed that if the cost functions are with some simple structures, one only needs to update the dual variables by some simple operations, which leads to the reduction of overall computational complexity.


\appendix{}

\subsection{Proof of Lemma \ref{lam2}}\label{lam2p}

By Lemma \ref{l1}, $ \nabla f^{\diamond}_i$ is Lipschitz continuous with Lipschitz constant $\frac{1}{\sigma_i}$, which means $ \| \nabla_{\mathbf{v}} f^{\diamond}_i(\mathbf{H}_i\mathbf{v}) - \nabla_{\mathbf{u}} f^{\diamond}_i(\mathbf{H}_i\mathbf{u}) \|
= \| \mathbf{H}_i^{\top} \nabla_{\mathbf{H}_i\mathbf{v}} f^{\diamond}_i(\mathbf{H}_i\mathbf{v}) - \mathbf{H}_i^{\top} \nabla_{\mathbf{H}_i\mathbf{u}} f^{\diamond}_i(\mathbf{H}_i\mathbf{u}) \| \leq \| \mathbf{H}_i \| \| \nabla_{\mathbf{H}_i\mathbf{v}} f^{\diamond}_i(\mathbf{H}_i\mathbf{v}) - \nabla_{\mathbf{H}_i\mathbf{u}} f^{\diamond}_i(\mathbf{H}_i\mathbf{u}) \|  \leq \frac{\| \mathbf{H}_i \| }{\sigma_i} \| \mathbf{H}_i \mathbf{v} - \mathbf{H}_i \mathbf{u} \| \leq  \frac{\| \mathbf{H}_i  \|^2 }{\sigma_i}\|  \mathbf{v}- \mathbf{u} \|
=  h_i \|  \mathbf{v}- \mathbf{u} \|$. Then, $ \nabla_{\bm{\lambda}_i} f^{\diamond}_i(\mathbf{H}_i\bm{\lambda}_i)$ is Lipschitz continuous with constant $h_i$ and, therefore, $\nabla_{\bm{\lambda}_i} p_i (\bm{\lambda}_i)=\nabla_{\bm{\lambda}_i} f^{\diamond}_i(\mathbf{H}_i\bm{\lambda}_i) + \kappa_i \mathbf{E}^{\top}$ is also Lipschitz continuous with constant $h_i$.


\subsection{Proof of Theorem \ref{la1}}\label{la1p}

By the first-order optimality condition of (\ref{f1}), we have
\begin{align}\label{e-1}
\mathbf{0}  \in & \partial_{\bm{\lambda}} Q (\bm{\lambda}^{t+1}) + c^{-1}(\bm{\lambda}^{t+1} - \bm{\lambda}^t) \nonumber \\
&  + \nabla_{\bm{\lambda}} P (\bm{\lambda}^t)  + \mathbf{M}^{\top} \bm{\xi}^t +  \gamma \mathbf{M}^{\top} \mathbf{M}  \bm{\lambda}^t \nonumber \\
 = &  \partial_{\bm{\lambda}} Q (\bm{\lambda}^{t+1}) - c^{-1}(\bm{\lambda}^t- \bm{\lambda}^{t+1}) + \nabla_{\bm{\lambda}} P (\bm{\lambda}^t) \nonumber \\
& +  \mathbf{M}^{\top} \bm{\xi}^{t+1} - \gamma \mathbf{M}^{\top} \mathbf{M}  \bm{\lambda}^{t+1}  +  \gamma \mathbf{M}^{\top} \mathbf{M}  \bm{\lambda}^t.
\end{align}
By the convexity of $Q(\bm{\lambda})$, we have
\begin{align}\label{e2}
Q  (\bm{\lambda}) & -  Q(\bm{\lambda}^{t+1}) \geq  c^{-1} (\bm{\lambda}-\bm{\lambda}^{t+1})^{\top} (\bm{\lambda}^t-\bm{\lambda}^{t+1}) \nonumber \\
& - (\bm{\lambda}-\bm{\lambda}^{t+1})^{\top} \nabla_{\bm{\lambda}} P(\bm{\lambda}^t) -(\bm{\lambda}-\bm{\lambda}^{t+1})^{\top} \mathbf{M}^{\top} \bm{\xi}^{t+1} \nonumber \\
&  + \gamma (\bm{\lambda}-\bm{\lambda}^{t+1})^{\top} \mathbf{M}^{\top} \mathbf{M}  (\bm{\lambda}^{t+1} -\bm{\lambda}^t ).
\end{align}
By the convexity and Lipschitz continuous differentiability of $p_i$, we have
\begin{align}\label{e3}
& (  \bm{\lambda} - \bm{\lambda}^{t+1})^{\top} \nabla_{\bm{\lambda}} P (\bm{\lambda}^t) \nonumber \\
= & \sum_{i \in {V}} (\bm{\lambda}_i - \bm{\lambda}_i^t)^{\top} \nabla_{\bm{\lambda}_i} p_i (\bm{\lambda}_i^t) + \sum_{i \in {V}} (\bm{\lambda}_i^t - \bm{\lambda}_i^{t+1})^{\top} \nabla_{\bm{\lambda}_i} p_i (\bm{\lambda}_i^t) \nonumber \\
\leq & \sum_{i \in {V}}( p_i(\bm{\lambda}_i) -  p_i(\bm{\lambda}_i^t) )  +  \sum_{i \in {V}}( p_i(\bm{\lambda}_i^t) -  p_i(\bm{\lambda}_i^{t+1}) ) \nonumber \\
& + \sum_{i \in {V}} \frac{h_i}{2}\| \bm{\lambda}_i^t - \bm{\lambda}_i^{t+1} \|^2 \nonumber  \\
\leq & P(\bm{\lambda})  -  P(\bm{\lambda}^{t+1}) + \frac{{h}}{2} \| \bm{\lambda}^t  - \bm{\lambda}^{t+1} \|^2.
\end{align}
By (\ref{f2}), we have $\mathbf{0}  =  \gamma^{-1}(\bm{\xi}^t-\bm{\xi}^{t+1}) + \mathbf{M} \bm{\lambda}^{t+1}$.
Multiplying the both sides of the above equality by $(\bm{\xi}-\bm{\xi}^{t+1})^{\top}$ gives
\begin{align}\label{r1-1}
 \gamma^{-1} (\bm{\xi} & -\bm{\xi}^{t+1})^{\top} (\bm{\xi}^t - \bm{\xi}^{t+1}) +  (\bm{\xi}-\bm{\xi}^{t+1})^{\top} \mathbf{M}\bm{\lambda}^{t+1} = 0.
\end{align}
By adding (\ref{e2}) and (\ref{e3}) together from the both sides, we have
\begin{align}\label{r1}
& \Phi (\bm{\lambda}^{t+1}) - \Phi(\bm{\lambda})
 \leq  - c^{-1} (\bm{\lambda} -\bm{\lambda}^{t+1})^{\top}  (\bm{\lambda}^t -\bm{\lambda}^{t+1})  \nonumber\\
& +(\bm{\lambda} -\bm{\lambda}^{t+1})^{\top} \mathbf{M}^{\top} \bm{\xi}^{t+1}  + \frac{{h}}{2} \| \bm{\lambda}^t - \bm{\lambda}^{t+1}\|^2 \nonumber\\
& - \gamma (\bm{\lambda} -\bm{\lambda}^{t+1})^{\top} \mathbf{M}^{\top} \mathbf{M}  (\bm{\lambda}^{t+1} -\bm{\lambda}^t ) \nonumber\\
= & - c^{-1} (\bm{\lambda} -\bm{\lambda}^{t+1})^{\top} (\bm{\lambda}^t -\bm{\lambda}^{t+1}) \nonumber\\
& - \gamma^{-1} (\bm{\xi} -\bm{\xi}^{t+1})^{\top} (\bm{\xi}^t -\bm{\xi}^{t+1}) - (\bm{\xi}-\bm{\xi}^{t+1})^{\top}\mathbf{M} \bm{\lambda}^{t+1}  \nonumber\\
& +\bm{\xi}^{t+1\top} \mathbf{M} \bm{\lambda} -\bm{\xi}^{t+1\top} \mathbf{M} \bm{\lambda}^{t+1} + \frac{{h}}{2} \| \bm{\lambda}^t - \bm{\lambda}^{t+1}\|^2 \nonumber \\
& - \gamma (\bm{\lambda} -\bm{\lambda}^{t+1})^{\top} \mathbf{M}^{\top} \mathbf{M}  (\bm{\lambda}^{t+1} -\bm{\lambda}^t )  \nonumber\\
= & \frac{1}{2c} (\|  \bm{\lambda} -\bm{\lambda}^t \|^2 - \|  \bm{\lambda} -\bm{\lambda}^{t+1} \|^2   -  \| \bm{\lambda}^t -\bm{\lambda}^{t+1}\|^2) \nonumber \\
& + \frac{1}{2\gamma} (\|  \bm{\xi} -\bm{\xi}^t \|^2 - \|  \bm{\xi} -\bm{\xi}^{t+1}\|^2 -  \| \bm{\xi}^t -\bm{\xi}^{t+1}\|^2 )  \nonumber\\
& +\bm{\xi}^{t+1\top} \mathbf{M} \bm{\lambda} - \bm{\xi}^{\top} \mathbf{M} \bm{\lambda}^{t+1} + \frac{{h}}{2} \| \bm{\lambda}^t - \bm{\lambda}^{t+1}\|^2 -  \|  \bm{\lambda} \nonumber\\
&  -\bm{\lambda}^t \|^2_{\frac{\gamma }{2} \mathbf{M}^{\top} \mathbf{M} } + \|  \bm{\lambda} -\bm{\lambda}^{t+1} \|^2_{\frac{\gamma }{2}\mathbf{M}^{\top} \mathbf{M} }  +  \| \bm{\lambda}^t -\bm{\lambda}^{t+1}\|^2_{\frac{\gamma }{2} \mathbf{M}^{\top} \mathbf{M} }  \nonumber \\
= & \|  \bm{\lambda} -\bm{\lambda}^t \|^2_{ \frac{1}{2c} \mathbf{I}_{NM+NB}- \frac{\gamma }{2}\mathbf{M}^{\top} \mathbf{M} } \nonumber\\
& - \|  \bm{\lambda} -\bm{\lambda}^{t+1} \|^2_{\frac{1}{2c} \mathbf{I}_{NM+NB}- \frac{\gamma }{2}\mathbf{M}^{\top} \mathbf{M} }  \nonumber\\
& -  \| \bm{\lambda}^t -\bm{\lambda}^{t+1} \|^2_{\left( \frac{1}{2c}- \frac{{h}}{2}\right) \mathbf{I}_{NM+NB} - \frac{\gamma }{2}\mathbf{M}^{\top} \mathbf{M} }  \nonumber\\
&  + \frac{1}{2\gamma} (\|  \bm{\xi} -\bm{\xi}^t\|^2 - \|  \bm{\xi} -\bm{\xi}^{t+1} \|^2 - \| \bm{\xi}^t  -\bm{\xi}^{t+1}\|^2 ) \nonumber\\
&  + \bm{\xi}^{t+1\top} \mathbf{M} \bm{\lambda} - \bm{\xi}^{\top} \mathbf{M} \bm{\lambda}^{t+1},
\end{align}
where we use (\ref{r1-1}) in the first equality and the second equality holds with $\mathbf{v}^{\top} \mathbf{u} = \frac{1}{2} (\| \mathbf{v} \|^2 + \| \mathbf{u} \|^2 - \| \mathbf{v} - \mathbf{u} \|^2) $.

Let $\bm{\xi} = \bm{\xi}^*$ and $\bm{\lambda} = \bm{\lambda}^*$ and rearrange (\ref{r1}), then we have
\begin{align}\label{r2}
& \Phi (\bm{\lambda}^{t+1}) - \Phi(\bm{\lambda}^*) +  \bm{\xi}^{*\top} \mathbf{M} \bm{\lambda}^{t+1} \nonumber\\
\leq & \|  \bm{\lambda}^* -\bm{\lambda}^t \|^2_{ \frac{1}{2c} \mathbf{I}_{NM+NB}- \frac{\gamma }{2}\mathbf{M}^{\top} \mathbf{M}} \nonumber\\
& - \|  \bm{\lambda}^* -\bm{\lambda}^{t+1} \|^2_{ \frac{1}{2c} \mathbf{I}_{NM+NB}- \frac{\gamma }{2}\mathbf{M}^{\top} \mathbf{M}} \nonumber\\
& -  \| \bm{\lambda}^t -\bm{\lambda}^{t+1}\|^2_{ \left( \frac{1}{2c}- \frac{{h}}{2}\right) \mathbf{I}_{NM+NB} - \frac{\gamma }{2}\mathbf{M}^{\top} \mathbf{M} } \nonumber\\
& + \frac{1}{2\gamma}(\|  \bm{\xi}^* -\bm{\xi}^t \|^2 - \|  \bm{\xi}^* -\bm{\xi}^{t+1} \|^2  - \| \bm{\xi}^t  - \bm{\xi}^{t+1} \|^2),
\end{align}
where (\ref{k2}) is used. By combining (\ref{fp1}), (\ref{sad}) and (\ref{k2}), we have
\begin{align}\label{pp1}
& \Phi( {\bm{\lambda}} ) - \Phi(\bm{\lambda}^*) +   \bm{\xi}^{*\top}\mathbf{M}  {\bm{\lambda}}  \geq 0.
\end{align}
Based on (\ref{r2}) and (\ref{pp1}), we have $a^{t+1} + b^t \leq a^t $, where $a^t = \|  \bm{\lambda}^* -\bm{\lambda}^t \|^2_{ \frac{1}{2c} \mathbf{I}_{NM+NB}- \frac{\gamma }{2}\mathbf{M}^{\top} \mathbf{M}}  + \frac{1}{2\gamma} \|  \bm{\xi}^* -\bm{\xi}^t \|^2,
b^t =  \| \bm{\lambda}^t -\bm{\lambda}^{t+1}\|^2_{ \left( \frac{1}{2c}- \frac{{h}}{2}\right) \mathbf{I}_{NM+NB} - \frac{\gamma }{2}\mathbf{M}^{\top} \mathbf{M}}  + \frac{1}{2\gamma} \| \bm{\xi}^t  - \bm{\xi}^{t+1} \|^2$. 
Based on (\ref{37}) and Thm. 3.1 in \cite{horn1998eigenvalue} (note that $\mathbf{I}_{NM+NB}$ and $\mathbf{M}^{\top} \mathbf{M}$ are Hermitian), we have $  0 \leq  c^{-1} - {h}  - \overline{\tau}(\gamma \mathbf{M}^{\top} \mathbf{M} )  <  c^{-1} - \overline{\tau}(\gamma \mathbf{M}^{\top} \mathbf{M} )  = \underline{\tau}\left(c^{-1}\mathbf{I}_{NM+NB} \right) - \overline{\tau}\left( \gamma \mathbf{M}^{\top} \mathbf{M} \right)   = \underline{\tau}\left(c^{-1}\mathbf{I}_{NM+NB} \right) + \underline{\tau}\left(-\gamma \mathbf{M}^{\top} \mathbf{M} \right) \leq \underline{\tau}\left(c^{-1}\mathbf{I}_{NM+NB}- \gamma \mathbf{M}^{\top} \mathbf{M} \right)$ and $0 \leq c^{-1} - {h} - \overline{\tau}(\gamma \mathbf{M}^{\top} \mathbf{M} ) = \underline{\tau}\left(( c^{-1}- {h}) \mathbf{I}_{NM+NB} \right) - \overline{\tau}\left( \gamma \mathbf{M}^{\top} \mathbf{M} \right) = \underline{\tau}\left(( c^{-1}- {h}) \mathbf{I}_{NM+NB} \right) + \underline{\tau}\left( - \gamma \mathbf{M}^{\top} \mathbf{M} \right)  \leq  \underline{\tau}\left(( c^{-1}- {h}) \mathbf{I}_{NM+NB} - \gamma \mathbf{M}^{\top} \mathbf{M} \right)$, which means
\begin{align}
& c^{-1}\mathbf{I}_{NM+NB} - \gamma \mathbf{M}^{\top} \mathbf{M}  \succ 0, \label{ee1} \\
&  ( c^{-1}- {h}) \mathbf{I}_{NM+NB} - \gamma \mathbf{M}^{\top} \mathbf{M}  \succeq 0. \label{ee2}
\end{align}
Then we can have $\lim_{t \rightarrow \infty} a^t$ exists and $\sum_{t=0}^{\infty}  b^t  < \infty$ \cite[Lemma 1]{chang2014distributed}. Then we can have the following results.
\begin{enumerate}
  \item Sequences $\{\bm{\lambda}^t\}_{t \in \mathbb{N}}$ and $\{\bm{\xi}^t\}_{t \in \mathbb{N}}$ are bounded. Then by Bolzano-Weierstrass Theorem \cite{bartle2000introduction}, there exists an increasing sequence $\{t_n\}_{n \in \mathbb{N}_+} \subseteq \mathbb{N}$ such that $\lim_{n \rightarrow \infty} \bm{\lambda}^{t_n} \triangleq \hat{\bm{\lambda}}$ and $\lim_{n \rightarrow \infty} \bm{\xi}^{t_n} \triangleq \hat{\bm{\xi}}$.
  \item $\lim_{t \rightarrow \infty} \| \bm{\lambda}^{t+1} - \bm{\lambda}^t \| = 0$ and $\lim_{t \rightarrow \infty} \| \bm{\xi}^{t+1} - \bm{\xi}^t \| = 0$, which means $\lim_{n \rightarrow \infty} \bm{\lambda}^{t_n+1} = \lim_{n \rightarrow \infty} \bm{\lambda}^{t_n-1} = \hat{\bm{\lambda}}$ and $\lim_{n \rightarrow \infty} \bm{\xi}^{t_n+1} = \lim_{n \rightarrow \infty} \bm{\xi}^{t_n-1} = \hat{\bm{\xi}}$.
\end{enumerate}

Therefore, $\mathbf{M}\hat{\bm{\lambda}} = \lim_{n \rightarrow \infty} \mathbf{M}\bm{\lambda}^{t_n+1} =  \lim_{n \rightarrow \infty} \gamma^{-1}\left( \bm{\xi}^{t_n+1}-\bm{\xi}^{t_n}\right) =  \mathbf{0}$, Then by Thm. 24.4 in \cite{rockafellar2015convex} and taking the subsequential limit of the both sides of (\ref{e-1}) along the instants $\{t_n\}_{n \in \mathbb{N}_+}$, we have $\mathbf{0} \in  \partial_{\bm{\lambda}} Q (\hat{\bm{\lambda}} ) + \nabla_{\bm{\lambda}} P (\hat{\bm{\lambda}}) + \mathbf{M}^{\top} \hat{\bm{\xi}}$. Therefore, $(\hat{\bm{\lambda}},\hat{\bm{\xi}})$ is a saddle point defined by (\ref{k1}) and (\ref{k2}). Define $a^{t_n} = \|  \hat{\bm{\lambda}} -\bm{\lambda}^{t_n} \|^2_{\frac{1}{2c} \mathbf{I}_{NM+NB} - \frac{\gamma}{2} \mathbf{M}^{\top} \mathbf{M} }  + \frac{1}{2\gamma} \|  \hat{\bm{\xi}} -\bm{\xi}^{t_n} \|^2$. Then $\lim_{n \rightarrow \infty} a^{t_n} = 0$. Since $\lim_{t \rightarrow \infty} a^t$ exists, then we have $\lim_{t \rightarrow \infty} a^t = \lim_{n \rightarrow \infty} a^{t_n} = 0$, which means $\lim_{t \rightarrow \infty} \bm{\lambda}^t = \hat{\bm{\lambda}}$ and $\lim_{t \rightarrow \infty} \bm{\xi}^t = \hat{\bm{\xi}}$. This proves the theorem.

\subsection{{Proof of Theorem \ref{th1}}}\label{th1p}

Note that (\ref{r1}) holds for all $\bm{\lambda}$ and $\bm{\xi}$. 
By letting $\bm{\lambda}=\bm{\lambda}^*$ and $ \bm{\xi} = 2\| \bm{\xi}^*\| \frac{\mathbf{M} \overline{\bm{\lambda}}^{T+1}}{\| \mathbf{M} \overline{\bm{\lambda}}^{T+1} \| }$ in (\ref{r1}), we have $\Phi (\bm{\lambda}^{t+1}) - \Phi(\bm{\lambda}^*) +  2\| \bm{\xi}^*\| \frac{(\mathbf{M} \overline{\bm{\lambda}}^{T+1})^{\top}}{\| \mathbf{M} \overline{\bm{\lambda}}^{T+1} \| } \mathbf{M} \bm{\lambda}^{t+1}
 \leq \|  \bm{\lambda}^* -\bm{\lambda}^t \|^2_{\frac{1}{2c} \mathbf{I}_{NM+NB}- \frac{\gamma}{2 }\mathbf{M}^{\top} \mathbf{M} }   -\|  \bm{\lambda}^* -\bm{\lambda}^{t+1} \|^2_{\frac{1}{2c} \mathbf{I}_{NM+NB} - \frac{\gamma}{2} \mathbf{M}^{\top} \mathbf{M} } + \frac{1}{2\gamma}  \left\|  2\| \bm{\xi}^*\| \frac{\mathbf{M} \overline{\bm{\lambda}}^{T+1}}{\| \mathbf{M} \overline{\bm{\lambda}}^{T+1} \| }  -\bm{\xi}^t \right\|^2 - \frac{1}{2\gamma}  \left\|  2\| \bm{\xi}^*\| \frac{\mathbf{M} \overline{\bm{\lambda}}^{T+1}}{\| \mathbf{M} \overline{\bm{\lambda}}^{T+1} \| }  -\bm{\xi}^{t+1} \right\|^2$, where (\ref{ee2}) is considered. Summing up above inequality over $t=0,1,...,T$ gives $(T +1)(\Phi(\overline{\bm{\lambda}}^{T+1}) - \Phi(\bm{\lambda}^*) +  2\| \bm{\xi}^* \| \|\mathbf{M} \overline{\bm{\lambda}}^{T+1} \|)  \leq  \sum_{t=0}^{T} ( \Phi({\bm{\lambda}}^{t+1}) - \Phi(\bm{\lambda}^*) +  2\| \bm{\xi}^* \| \|\mathbf{M} \overline{\bm{\lambda}}^{T+1} \| ) \leq  \|  \bm{\lambda}^* -\bm{\lambda}^0 \|^2_{\frac{1}{2c} \mathbf{I}_{NM+NB} - \frac{\gamma}{2} \mathbf{M}^{\top} \mathbf{M} }  + \frac{1}{2\gamma} \left\|  2 \| \bm{\xi}^*\| \frac{\mathbf{M} \overline{\bm{\lambda}}^{T+1}}{\| \mathbf{M} \overline{\bm{\lambda}}^{T+1} \| } -\bm{\xi}^0 \right\|^2
 \leq  \|  \bm{\lambda}^* -\bm{\lambda}^0 \|^2_{\frac{1}{2c}\mathbf{I}_{NM+NB} - \frac{\gamma}{2} \mathbf{M}^{\top} \mathbf{M} } + \frac{4}{\gamma} \|  \bm{\xi}^* \|^2 + \frac{1}{\gamma} \| \bm{\xi}^0 \|^2
 \triangleq  \Theta$, where the first inequality is from the convexity of $\Phi$ and the third inequality is from Cauchy-Schwarz inequality. Therefore,
\begin{align}\label{sls1}
\Phi (\overline{\bm{\lambda}}^{T+1}) - \Phi(\bm{\lambda}^*)
\leq \frac{\Theta}{T+1}  - 2\| \bm{\xi}^* \| \|\mathbf{M} \overline{\bm{\lambda}}^{T+1} \| \leq \frac{\Theta}{T+1}.
\end{align}
By letting $\bm{\lambda}=\overline{\bm{\lambda}}^{T+1}$ in (\ref{pp1}), we have
\begin{align}\label{sls2}
\Phi(\overline{\bm{\lambda}}^{T+1}) - \Phi(\bm{\lambda}^*)\geq - \| \bm{\xi}^* \| \|\mathbf{M} \overline{\bm{\lambda}}^{T+1} \|.
\end{align}
By combining (\ref{sls1}) and (\ref{sls2}), we have
\begin{align}\label{sls2+1}
\| \bm{\xi}^* \| \| \mathbf{M} \overline{\bm{\lambda}}^{T+1} \| \leq  \frac{\Theta}{T+1}.
\end{align}
Then with (\ref{sls2}), we can further have
\begin{align}\label{sls3}
 \Phi  (\overline{\bm{\lambda}} ^{T+1})  - \Phi(\bm{\lambda}^*)  \geq  - \frac{\Theta}{T+1}.
\end{align}
The proof is completed based on (\ref{sls1}), (\ref{sls2+1}), and (\ref{sls3}).


\bibliographystyle{IEEEtran}
\small
\bibliography{1myref}

\end{document}